\documentclass[aos,MSNbibl,seceqn,nameyear,dvips]{arximspdf}
\usepackage{graphicx}

\doi{10.1214/13-AOS1137} %kopijuoti is PTS
\volume{41}
\issue{1}
\pubyear{2013}
\firstpage{1892}
\lastpage{1921}

\makeatletter
\newtheorem{theorem}{Theorem}[section] %Emily - need to include
\newtheorem{corollary}{Corollary}[section] %E %Emily - need to include
\newproclaim{remark}{Remark} %Emily - need to include
\def\cal{\mathcal}
\def\a{\alpha}
\def\B{\textrm{B}}
\def\be{\beta}
\def\bep{{\bar\ep}}
\def\cB{{\cal B}}
\def\cE{{\cal E}}
\def\cR{{\cal R}}
\def\cS{{\cal S}}
\def\cX{{\cal X}}
\def\cZ{{\cal Z}}
\def\cov{\operatorname{cov}}
\def\De{\Delta}
\def\ep{\varepsilon}
\def\excep{_{\mathrm{excep}}}
\def\ga{\gamma}
\def\hal{{\hat\a}}
\def\hbe{{\hat\be}}
\def\hd{{\hat d}}
\def\hep{{\hat\ep}}
\def\hf{{\hat f}}
\def\hg{{\hat g}}
\def\hpi{{\hat\pi}}
\def\hsi{{\hat\si}}
\def\hz{{\hat z}}
\def\IR{\mathbb{R}}
\def\ka{\kappa}
\def\la{\lambda}
\def\med{\operatorname{med}}
\def\part{\partial}
\def\quant{\mathrm{quant}}
\def\ra{\to}
\def\rai{\ra\infty}
\def\sgn{\operatorname{sgn}}
\def\si{\sigma}
\def\slim{_{\mathrm{lim}}}
\def\tep{{\tilde\ep}}
\def\var{\operatorname{var}}
\makeatother

\begin{document}
\begin{frontmatter}

\title{A simple bootstrap method for constructing nonparametric
confidence bands for functions}
\runtitle{Confidence bands}

\begin{aug}
\author[a]{\fnms{Peter} \snm{Hall}\corref{}\thanksref{t1}\ead[label=PH]{halpstat@ms.unimelb.edu.au}}
\and
\author[b]{\fnms{Joel} \snm{Horowitz}\thanksref{t2}\ead[label=JH]{joel-horowitz@northwestern.edu}}
\affiliation{University of Melbourne and University of California, and~Northwestern~University}

\thankstext{t1}{Supported by ARC and NSF grants.}
\thankstext{t2}{Supported by NSF Grant SES-0817552.}
\address[a]{Department of Mathematics and Statistics\\
University of Melbourne \\
VIC 3010\\
Australia\\
and\\
Department of Statistics\\
University of California \\
Davis, California 95616\\
USA\\
\printead{PH}}

\address[b]{Department of Economics \\
Northwestern University\\
2001 Sheridan Road\\
Evanston, Illinois 60208\\
USA\\
\printead{JH}}

\runauthor{P. Hall and J. Horowitz}
\end{aug}

% HISTORY:
\received{\smonth{1} \syear{2013}}
\revised{\smonth{6} \syear{2013}}

% ABSTRACT
%
\begin{abstract}
Standard approaches to constructing nonparametric confidence bands for
functions are frustrated by the impact of bias, which generally is not
estimated consistently when using the bootstrap and conventionally
smoothed function estimators. To overcome this problem it is common
practice to either undersmooth, so as to reduce the impact of bias, or
oversmooth, and thereby introduce an explicit or implicit bias
estimator. However, these approaches, and others based on nonstandard
smoothing methods, complicate the process of inference, for example, by
requiring the choice of new, unconventional smoothing parameters and,
in the case of undersmoothing, producing relatively wide bands. In this
paper we suggest a new approach, which exploits to our advantage one of
the difficulties that, in the past, has prevented an attractive
solution to the problem---the fact that the standard bootstrap bias
estimator suffers from relatively high-frequency stochastic error. The
high frequency, together with a technique based on quantiles, can be
exploited to dampen down the stochastic error term, leading to
relatively narrow, simple-to-construct confidence bands.
\end{abstract}

% KEYWORDS

\end{frontmatter}

%s1 #&#
\section{Introduction}\label{Intro}

%s1.1 #&#
\subsection{Motivation}\label{Motivation}
There is an extensive literature, summarised in Section~\ref{LitReview}
below, on constructing nonparametric confidence bands for functions.
However, this work generally does not suggest practical solutions to
the critical problem of choosing tuning parameters, for example,
smoothing parameters or the nominal coverage level of the confidence
band, to ensure a high degree of coverage accuracy or to produce bands
that err on the side of conservatism. In this paper we suggest new,
simple bootstrap methods for constructing confidence bands using
conventional smoothing parameter choices.

In particular, our approach does not require a nonstandard smoothing
parameter. The basic algorithm requires only a single application of
the bootstrap, although a more refined, double bootstrap technique is
also suggested. The greater part of our attention is directed to
regression problems, but we also discuss the application of our methods
to constructing confidence bands for density functions.

The resulting confidence regions depend on choice of two parameters $\a
$ and $\xi$, in the range $0<\a,\xi<1$, and the methodology results in
confidence bands that, asymptotically, cover the regression mean at $x$
with probability at least $1-\a$, for at least a proportion $1-\xi$ of
values of $x$. In particular, the bands are pointwise, rather than
simultaneous. Pointwise bands are more popular with practitioners and
are the subject of a substantial majority of research on nonparametric
confidence bands for functions.

%s1.2 #&#
\subsection{Features of our approach, and competing methods}\label{Features}
The ``exceptional'' $100 \xi\%$ of points that are not covered are
typically close to the locations of peaks and troughs, and so are
discernible from a simple estimate of the regression mean. Their
location can also be determined using a theoretical analysis---points
near peaks and troughs potentially cause difficulties because of bias.
See Section~\ref{Values} for theoretical details, and Section~\ref{NumProps} for numerical examples.

Our approach accommodates bias by increasing the width of confidence
bands. However, the amount by which we increase width is no greater
than a constant factor, rather than the polynomial amount (as a
function of $n$) associated with most suggestions for undersmoothing.

Methods based on either under- or oversmoothing are recommended often
in the literature. However, there are no empirical techniques, where
the data determine the amount of smoothing, that are used even
moderately widely in either case. In particular, although theoretical
arguments demonstrate clearly the advantages of under- or oversmoothing
if appropriate smoothing parameters are chosen, there are no
attractive, effective empirical ways of selecting those quantities.
Indeed, it is not uncommon to suggest that the issue be avoided
altogether, by ignoring the effects of bias. For example, this approach
is recommended in textbooks; see \citet{r53}, pages 133ff, who refer to
the resulting bands as ``variability bands,'' and \citet{r12}, pages 79--80, who suggest plotting many realisations of
bootstrapped curve estimators without bias corrections.

In addition to needing unavailable bandwidth choice methods, the
drawbacks of undersmoothing include the fact that the confidence bands
become both wider and more wiggly as the amount of undersmoothing
increases. The increase in wiggliness is so great that, unless sample
size is very large, the coverage accuracy does not necessarily improve
as the amount of undersmoothing increases. Details are given in
Section~\ref{NumProps}.

Wiggliness can likewise be a problem for bands that result from using
oversmoothing to remove bias explicitly. Here the relatively high level
of variability from which function derivative estimators suffer means
that the confidence bands may again oscillate significantly, and can be
difficult to interpret. These results, and those reported in the
previous paragraph, are for optimal choices of the amount of under- or
oversmoothing. In practice the amount has to be chosen empirically, and
that introduces additional noise, which further reduces performance.

%s1.3 #&#
\subsection{Intuition}\label{Intuition}
Our methodology exploits, to our advantage, a difficulty that in the
past has hindered a simple solution to the confidence band problem. To
explain how, we note first that if nonparametric function estimators
are constructed in a conventional manner then their bias is of the same
order as their error about the mean, and accommodating the bias has
been a major obstacle to achieving good coverage accuracy. Various
methods, based on conventional smoothing parameters, can be used to
estimate the bias and reduce its impact, but the bias estimators fail
to be consistent, not least because the stochastic noise from which
they suffer is highly erratic. (In the case of kernel methods, the
frequency of the noise is proportional to the inverse of the
bandwidth.) However, as we show in this paper, this erratic behaviour
is actually advantageous, since if we average over it, then we can
largely eliminate the negative impact that it has on the bias
estimation problem. We do the averaging implicitly, not by computing
means but by working with quantiles of the ``distribution'' of coverage.

%s1.4 #&#
\subsection{Literature review}\label{LitReview}
We shall summarise previous work largely in terms of whether it
involved undersmoothing or oversmoothing; the technique suggested in
the present paper is almost unique in that it requires neither of these
approaches. \citet{r28}, \citet{r31}, \citet{r21}, \citet{r13},
\citet
{r56}, \citet{r29} and \citet{r60} suggested methods based on
oversmoothing, using either implicit or explicit bias correction.
\citet
{r27} also used explicit bias correction, in the sense that their bands
required a known bound on an appropriate derivative of the target
function. \citet{r3}, \citet{r22}, \citet{r26}, \citet{r47}, \citet
{r8}, \citet{r48}, \citet{r49}, \citet{r9} (in the context of
hypothesis testing), \citet{r10}, \citet{r30} and \citet{r38} employed
methods that involve undersmoothing. There is also a theoretical
literature which addresses the bias issue through consideration of the
technical function class from which a regression mean or density came;
see, for example, \citet{r37} and \citet{r18}. This work sometimes
involves confidence balls, rather than bands, and in that respect is
connected to research such as that of \citet{r14} and \citet{r17}.
\citet{r59} considered spline and Bayesian methods. The notion of
``honest'' confidence bands, which have guaranteed coverage for a rich
class of functions, was pioneered by \citet{r35}. Recent contributions
include those of \citet{r7}, \citet{r19} and \citet{r32}.

%s2 #&#
\section{Methodology}\label{Methodology}

%s2.1 #&#
\subsection{Model}\label{Model}
Suppose we observe data pairs in a sample $\cZ=\{(X_i,Y_i), 1\leq
i\leq n\}$, generated by the model
%
%
%e2.1 #&#
%
\begin{equation}
Y_i=g(X_i)+\ep_i,\label{eq:2.1}
\end{equation}
where the experimental errors $\ep_i$ are independent and identically
distributed with finite variance and zero mean conditional on $X$. Our
aim is to construct a pointwise confidence band for the true $g$ in a
closed, bounded region $\cR$. A more elaborate, heteroscedastic model
will be discussed in Section~\ref{3Remarks}; we omit it here only for
the sake of simplicity. We interpret $g(x)$ in the conventional
regression manner, as $E(Y | X=x)$, but our theoretical analysis takes
account of the fact that although we condition on the $X_i$s at this
point we consider that they originated as random variables, with
density $f_X$.

%s2.2 #&#
\subsection{Properties of function estimators and conventional
confidence bands}\label{Properties}
Let~$\hg$ denote a conventional estimator of $g$. We assume that $\hg$
incorporates smoothing parameters computed empirically from the data,
using for example cross-validation or a plug-in rule, and that the
variance of $\hg$ can be estimated consistently by $s(\cX)^2 \hsi^2$,
where $s(\cX)$ is a known function of the set of design points $\cX=\{
X_1,\ldots,X_n\}$ and the smoothing parameters, and $\hsi^2$ is an
estimator of the variance, $\si^2$, of the experimental errors $\ep_i$,
computed from the dataset $\cZ$. The case of heteroscedasticity is
readily accommodated too; see Section~\ref{3Remarks}. We write $\hg
^*$
for the version of $\hg$ computed using a conventional bootstrap
argument. For details of the construction of $\hg^*$, see step~4 of
the algorithm in Section~\ref{Algorithm}.

The smoothing parameters used for $\hg$ would generally be chosen to
optimise a measure of accuracy, for example, in a weighted $L_p$ metric
where $1\leq p<\infty$, and we shall make this assumption implicitly in
the discussion below. In particular, it implies that the asymptotic
effect of bias, for example, as represented by the term $b(x)$ in (\ref
{eq:2.4}) below, is finite and typically nonzero.

An asymptotic, symmetric confidence band for $g$, constructed naively
without considering bias, and with nominal coverage $1-\a$, has the form
%
%
%e2.2 #&#
%
\begin{eqnarray}\label{eq:2.2}
&&\cB(\a)= \bigl\{(x,y)\dvtx x\in\cR, \hg(x)-s(\cX) (x) \hsi z_{1-(\a/2)}\leq
y
\nonumber
\\[-8pt]
\\[-8pt]
\nonumber
&&\hspace*{90pt}\qquad\leq\hg(x)+s(\cX) (x) \hsi z_{1-(\a/2)} \bigr\},
\end{eqnarray}
where $z_\be=\Phi^{-1}(\be)$ is the $\be$-level critical point of the
standard normal distribution, and $\Phi$ is the standard normal
distribution function. Unfortunately, the coverage of $\cB(\a)$ at a
point $x$, given by
%
%
%e2.3 #&#
%
\begin{equation}
\pi(x,\a)=P \bigl\{ \bigl(x,g(x) \bigr)\in\cB(\a) \bigr\},\label{eq:2.3}
\end{equation}
is usually incorrect even in an asymptotic sense, and in fact the band
typically undercovers, often seriously, in the limit as $n\rai$. The
reason is that the bias of $\hg$, as an estimator of $g$, is of the
same size as the estimator's stochastic error, and the confidence band
allows only for the latter type of error. As a result the limit, as
$n\rai$, of the coverage of the band is given by
%
%
%e2.4 #&#
%
\begin{equation}
\pi\slim(x,\a)=\lim_{n\rai} \pi(x,\a) =\Phi \bigl\{z+b(x) \bigr
\}- \Phi \bigl\{-z+b(x) \bigr\},\label{eq:2.4}
\end{equation}
where $z=z_{1-(\a/2)}$ and $b(x)$ describes the asymptotic effect that
bias has on coverage. [A formula for $b(x)$ in a general multivariate
setting is given in (\ref{eq:4.7}), and a formula in the univariate
case is provided in Section~\ref{Values}.] The right-hand side of~(\ref
{eq:2.4}) equals $\Phi(z)-\Phi(-z)=1-\a$ if and only if $b(x)=0$. For
all other values of $b(x)$, $\pi\slim(x,\a)<1-\a$. This explains why
the band at (\ref{eq:2.2}) almost always undercovers unless some sort
of bias correction is used.

The band potentially can be recalibrated, using the bootstrap, to
correct for coverage errors caused by bias, but now another issue
causes difficulty: the standard bootstrap estimator of bias, $E\{\hg
^*
(x) | \cZ\}-\hg(x)$, is inconsistent, in the sense that the ratio of
the estimated bias to its true value does not converge to 1 in
probability as $n\rai$. This time the problem is caused by the
stochastic error of the bias estimator; it is of the same size as the
bias itself. The problem can be addressed using an appropriately
oversmoothed version of $\hg$ when estimating bias, either explicitly
or implicitly, but the degree of oversmoothing has to be determined
from the data, and in practice this issue is awkward to resolve.
Alternatively, the estimator $\hg$ can be undersmoothed, so that the
influence of bias is reduced, but now the amount of undersmoothing has
to be determined, and that too is difficult. Moreover, confidence bands
computed from an appropriately undersmoothed $\hg$ are an order of
magnitude wider than those at (\ref{eq:2.2}), and so the undersmoothing
approach, although more popular than oversmoothing, is unattractive for
at least two reasons.

A simpler bootstrap technique, described in detail in the next section,
overcomes these problems.

%s2.3 #&#
\subsection{The algorithm}\label{Algorithm}

\subsubsection*{Step~1. Estimators of $g$ and $\si^2$}
Construct a conventional nonparametric estimator $\hg$ of $g$. Use a
standard empirical method (e.g., cross-validation or a plug-in rule),
designed to minimise mean $L_p$ error for some $p$ in the range $1\leq
p<\infty$, to choose the smoothing parameters on which $\hg$ depends.
For example, if the design is univariate then a local linear estimator
of $g(x)$ is given by
%
%
%e2.5 #&#
%
\begin{equation}
\hg(x)=\frac{1}{n}\sum_{i=1}^n
A_i(x) Y_i,\label{eq:2.5}
\end{equation}
where
%
%
%e2.6 #&#
%
\begin{equation}
A_i(x)=\frac{S_2(x)-\{(x-X_i)/h\} S_1(x)}{
S_0(x) S_2(x)-S_1(x)^2} K_i(x),\label{eq:2.6}
\end{equation}
$S_k(x)=n^{-1}\sum_i \{(x-X_i)/h\}^k K_i(x)$, $K_i(x)=h^{-1}K\{
(x-X_i)/h\}$, $K$ is a kernel function and $h$ is a bandwidth.

There is an extensive literature on computing estimators $\hsi^2$ of
the error variance $\si^2=\var(\ep)$; see, for example, \citet{r50},
\citet{r5}, \citet{r16}, M{\"u}ller and Stadtm{\"u}ller (\citeyear{r41,r42}), \citet{r24}, \citet{r25},
\citet{r55}, \citet{r46}, \citet{r43}, \citet{r11}, \citet{r15},
\citet
{r44}, \citet{r45}, \citet{r57}, \citet{r4}, \citet{r6}, and \citet
{r40}. It includes residual-based estimators, which we introduce at
(\ref{eq:2.8}) below, and methods based on differences and generalised
differences. An example of the latter approach, in the case of
univariate design, is the following estimator due to \citet{r50}:
%
%
%e2.7 #&#
%
\begin{equation}
\hsi^2=\frac{1}{2 (n-1)} \sum_{i=2}^n
(Y_{[i]}-Y_{[i-1]})^2,\label{eq:2.7}
\end{equation}
where $Y_{[i]}$ is the concomitant of $X_{(i)}$ and $X_{(1)}\leq\cdots
\leq X_{(n)}$ is the sequence of order statistics derived from the
design variables.

As in Section~\ref{Properties}, let $s(\cX)(x)^2 \hsi^2$ denote an
estimator of the variance of $\hg(x)$, where $s(\cX)(x)$ depends on the
data only through the design points, and $\hsi^2$ estimates error
variance, for example, being defined as at (\ref{eq:2.7}) or (\ref
{eq:2.8}). In the local linear example, introduced at (\ref{eq:2.5})
and (\ref{eq:2.6}), we take $s(\cX)(x)^2=\ka/\{nh \hf_X(x)\}$, where
$\ka=\int K^2$ and $\hf_X(x)=(nh_1)^{-1}\sum_{1\leq i\leq n} K_1\{
(x-X_i)/h_1\}$ is a standard kernel density estimator, potentially
constructed using a bandwidth $h_1$ and kernel $K_1$ different from
those used for $\hg$. There are many effective, empirical ways of
choosing $h_1$, and any of those can be used.

\subsubsection*{Step~2. Computing residuals}
Using the estimator $\hg$ from step (1), calculate initial residuals
$\tep_i=Y_i-\hg(X_i)$, put $\bep=n^{-1}\sum_i \tep_i$, and define the
centred residuals by $\hep_i=\tep_i-\bep$.

A conventional, residual-based estimator of $\si^2$, alternative to the
estimator at~(\ref{eq:2.7}), is
%
%
%e2.8 #&#
%
\begin{equation}
\hsi^2=\frac{1}{n}\sum_{i=1}^n
\hep_i^2.\label{eq:2.8}
\end{equation}
The estimator at (\ref{eq:2.7}) is root-$n$ consistent for $\si^2$,
whereas the estimator at (\ref{eq:2.8}) converges at a slower rate
unless an undersmoothed estimator of $\hg$ is used when computing the
residuals. This issue is immaterial to the theory in Section~\ref{TheoProps}, although it tends to make the estimator at (\ref{eq:2.7})
a little more attractive.

\subsubsection*{Step~3. Computing bootstrap resample}
Construct a resample $\cZ^*=\break\{(X_i,Y_i^*)$, $1\leq i\leq n\}$, where
$Y_i^*=\hg(X_i)+\ep_i^*$ and the $\ep_i^*$s are obtained by sampling
from $\hep_1,\ldots,\hep_n$ randomly, with replacement, conditional
on $\cX$. Note that, since regression is conventionally undertaken
conditional on the design sequence, then the $X_i$s are not resampled,
only the $Y_i$s.

\subsubsection*{Step~4. Bootstrap versions of $\hg$, $\hsi^2$ and
$\cB
(\a)$}
From the resample drawn in step~3, but using the same smoothing
parameter employed to construct~$\hg$, compute the bootstrap version
$\hg^*$ of $\hg$. (See Section~\ref{3Remarks} for discussion of the
smoothing parameter issue.) Let $\hsi^*{}^2$ denote the bootstrap
version of $\hsi^2$, obtained when the latter is computed from $\cZ
^*$
rather than $\cZ$, and construct the bootstrap version of $\cB(\a)$,
at (\ref{eq:2.2}),
%
%
%e2.9 #&#
%
\begin{eqnarray}\label{eq:2.9}
&&\cB^*(\a)= \bigl\{(x,y)\dvtx x\in\cR, \hg^*(x)-s(\cX) (x) \hsi^*
z_{1-(\a/2)}\leq y
\nonumber
\\[-8pt]
\\[-8pt]
\nonumber
&&\hspace*{95pt}\qquad\leq\hg^*(x)+s(\cX) (x) \hsi^*z_{1-(\a/2)} \bigr\}.
\end{eqnarray}
Note that $s(\cX)$ is exactly the same as in (\ref{eq:2.2}); again this
is a consequence of the fact that we are conducting inference
conditional on the design points.

If, as in the illustration in step~1, the design is univariate and
local linear estimators are employed, then $\hg^*(x)=n^{-1}\sum_{1\leq
i\leq n} A_i(x) Y_i^*$ where $A_i(x)$ is as at (\ref{eq:2.6}). The
bootstrap analogue of the variance formula (\ref{eq:2.7}) is $\hsi^*
{}^2=\{2 (n-1)\}^{-1}\sum_{2\leq i\leq n} (Y_{[i]}^*-Y_{[i-1]}^*
)^2$, where, if the $i$th largest order statistic $X_{(i)}$ equals
$X_j$, then $Y_{[i]}^*=\hg(X_j)+\ep_j^*$.

\subsubsection*{Step~5. Estimator of coverage error}
The bootstrap estimator $\hpi(x,\a)$ of the probability $\pi(x,\a)$
that $\cB(\a)$ covers $(x,g(x))$ is defined by
%
%
%e2.10 #&#
%
\begin{equation}
\hpi(x,\a)=P \bigl\{ \bigl(x,\hg(x) \bigr)\in\cB^*(\a) | \cX \bigr\},
\label{eq:2.10}
\end{equation}
and is computed, by Monte Carlo simulation, in the form
%
%
%e2.11 #&#
%
\begin{equation}
\frac{1}{B} \sum_{b=1}^B I \bigl\{
\bigl(x,\hg(x) \bigr)\in\cB_b^*(\a) \bigr\},\label{eq:2.11}
\end{equation}
where $I(\cE)$ denotes the indicator function of an event $\cE$, and
$\cB_b^*(\a)$ is the $b$th out of $B$ bootstrap replicates of $\cB
^*
(\a)$, where the latter is as at~(\ref{eq:2.9}). The estimator at~(\ref
{eq:2.10}) is completely conventional, and in particular, no additional
or nonstandard smoothing is needed.

\subsubsection*{Step~6. Constructing final confidence band}
Define $\hbe(x,\a_0)$ to be the solution, in $\a$, of $\hpi(x,\a
)=1-\a
_0$, and let $\hal_\xi(\a_0)$ denote the $\xi$-level quantile of points
in the set $\{\hbe(x,\a_0)\dvtx x\in\cR\}$. Specifically:
%
%
%e2.12 #&#
%
\begin{equation}
\begin{tabular}{p{280pt}} {take $\cR$ to be a subset of $
\IR^r$, superimpose on $\cR$ a regular, $r$-dimensional, rectangular
grid with edge width $\delta$, let $x_1,\ldots,x_N\in
\cR$ be the grid centres, let $\hal_\xi(\a_0,\delta)$
denote the $\xi$-level empirical quantile of the points $\hal(x_1,
\a_0),\ldots,\hal(x_N,\a_0)$, and, for $\xi
\in(0,1)$, let $\hal_\xi(\a_0)$ denote the limit infimum,
as $\delta\ra0$, of the sequence $\hal_\xi(\a_0,\delta
)$.} \end{tabular} %
\label{eq:2.12}
\end{equation}
(We use the limit infimum to avoid ambiguity, although under mild
conditions the limit exists.) For a value $\xi\in(0,{\frac{1}{2}}]$,
construct
the band $\cB\{\hal_\xi(\a_0)\}$. In practice we have found that taking
$1-\xi=0.9$ generally gives a slight to moderate degree of
conservatism, except for the exceptional points $x$ that comprise
asymptotically a fraction $\xi$ of $\cR$. Taking $1-\xi=0.95$ may be
warranted in the case of large samples.

%s2.4 #&#
\subsection{Three remarks on the algorithm}\label{3Remarks}
\begin{remark}[(Calibration)]\label{rem1}
In view of the undercoverage property discussed below (\ref{eq:2.4}),
we expect $\hbe(x,\a_0)$, defined in step~6, to be less than $\a_0$.
Equivalently, we anticipate that the nominal coverage of the band has
to be increased above $1-\a_0$ in order for the band to cover
$(x,g(x))$ with probability at least $1-\a_0$. Conventionally we would
employ $\hbe(x,\a_0)$ as the nominal level, but, owing to the large
amount of stochastic error in the bootstrap bias estimator that is used
implicitly in this technique, it produces confidence bands with poor
coverage accuracy. This motivates coverage correction by calibration,
along lines suggested by \citet{r20}, \citet{r1} and \citet{r36}, and
resulting in our use of the adjusted nominal level $\hal_\xi(\a_0)$,
defined in step~6.
\end{remark}

\begin{remark}[(Smoothing parameter for $\hg^*$)]\label{rem2}
An important aspect of step~4 is that we use the same empirical
smoothing parameters for both $\hg^*$ and $\hg$, even though,\vadjust{\goodbreak} in some
respects, it might seem appropriate to use a bootstrap version of the
smoothing parameters for $\hg$ when estimating $\hg^*$. However, since
smoothing parameters should be chosen to effect an optimal tradeoff
between bias and stochastic error, and the bias of $\hg$ is not
estimated accurately by the conventional bootstrap used in step~3
above, then the bootstrap versions of smoothing parameters, used to
construct $\hg^*$, are generally not asymptotically equivalent to
their counterparts used for $\hg$. This can cause difficulties. The
innate conservatism of our methodology accommodates the slightly
nonstandard smoothing parameter choice in step~4. Moreover, by not
having to recompute the bandwidth at every bootstrap step, we
substantially reduce computational labour.
\end{remark}

\begin{remark}[(Heteroscedasticity)]\label{rem3}
A heteroscedastic generalisation of the model at (\ref{eq:2.1}) has
the form
%
%
%e2.13 #&#
%
\begin{equation}
Y_i=g(X_i)+\si(X_i) \ep_i,
\label{eq:2.13}
\end{equation}
where the $\ep_i$s have zero mean and unit variance, and $\si(x)$ is a
nonnegative function that is estimated consistently by $\hsi(x)$, say,
computed from the dataset $\cZ$ using either parametric or
nonparametric methods. In this setting the variance of $\hg(x)$
generally can be estimated by $s(\cX)^2 \hsi(x)^2$, where $s(\cX)$ is
a known function of the design points, and the confidence band at (\ref
{eq:2.2}) should be replaced by
\begin{eqnarray*}
&&\cB(\a)= \bigl\{(x,y)\dvtx x\in\cR, \hg(x)-s(\cX) (x) \hsi(x) z_{1-(\a/2)}
\leq y \\
&&\hspace*{90pt}\qquad\leq\hg(x)+s(\cX) (x) \hsi(x) z_{1-(\a/2)} \bigr\}.
\end{eqnarray*}
The model for generating bootstrap data now has the form $Y_i^*=\hg
(X_i)+\hsi(X_i) \ep_i^*$, instead of $Y_i^*=\hg(X_i)+\ep_i^*$ in
step~4; and the $\ep_i^*$s are resampled conventionally from residual
approximations to the $\ep_i$s.

With these modifications, the algorithm described in steps 1--6 can be
implemented as before, and the resulting confidence bands have similar
properties. In particular, if we redefine $\cB^*(\a)$ by
\begin{eqnarray*}
&&\cB^*(\a)= \bigl\{(x,y)\dvtx x\in\cR, \hg^*(x)-s(\cX) (x) \hsi^*(x)
z_{1-(\a/2)}\leq y \\
&&\hspace*{95pt}\qquad\leq\hg^*(x)+s(\cX) (x) \hsi^*(x) z_{1-(\a/2)} \bigr\}
\end{eqnarray*}
[compare (\ref{eq:2.9})], and, using this new definition, continue to
define $\hpi(x,\a)$ as at~(\ref{eq:2.10}) [computed as at (\ref
{eq:2.11})]; and if we continue to define $\be=\hbe(x,\a_0)$ to be the
solution of $\hpi(x,\be)=1-\a_0$, and to define $\hal_\xi(\a_0)$
as in
(\ref{eq:2.12}); then the confidence band $\cB\{\hal_\xi(\a_0)\}$ is
asymptotically conservative for at least a proportion $1-\xi$ of values
$x\in\cR$. This approach can be justified intuitively as in
Appendix~B.1 in the supplementary file, noting that, in the context of
the model at (\ref{eq:2.13}), the expansion at (\B.1) in the supplement
should be replaced by
\begin{eqnarray*}
E \bigl\{\hg^*(x) | \cZ \bigr\}-\hg(x)&=&c_1 g''(x)
h^2+(nh)^{-1/2}\si(x) f_X(x)^{-1/2}W(x/h)
\\
&&{}+\mbox{negligible terms}.
\end{eqnarray*}
\end{remark}

%s2.5 #&#
\subsection{Percentile bootstrap confidence bands}\label{Percentile}
The methods discussed above are based on the symmetric, asymptotic
confidence band $\cB(\a)$, which in turn is founded on a normal
approximation. This approach is attractive because it requires only a
single application of the bootstrap for calibration, but it is
restrictive in that it dictates a conventional, symmetric ``template''
for the bands, because the normal model is symmetric. However,
particularly if we would prefer the bands to be placed asymmetrically
on either side of the estimator $\hg$ so as to reflect skewness of the
distribution of experimental errors, the initial confidence band $\cB
(\a
)$, at (\ref{eq:2.2}), can be constructed using bootstrap methods, and
a second iteration of the bootstrap, resulting in a double bootstrap
method, can be used to refine coverage accuracy. This allows us to use,
for example, equal-tailed intervals (where the amount of probability in
either tail is taken to be the same) and so-called ``shortest''
intervals (where the confidence interval is chosen to be as short as
possible, subject to having the desired nominal coverage). Of course,
one-sided intervals can be constructed using either a normal
approximation or a bootstrap approach, and our method carries over
without difficulty to those settings.

The first bootstrap implementation is undertaken using step~4 of the
algorithm in Section~\ref{Algorithm}, and allows us to define the
critical point $\hz_\be(x)$ by
%
%
%e2.14 #&#
%
\begin{equation}
P \bigl\{\hg^*(x)-\hg(x)\leq s(\cX) \hz_\be | \cZ \bigr\}=\be \label
{eq:2.14}
\end{equation}
for $\be\in(0,1)$. The confidence band $\cB(\a)$ is now re-defined as
%
%
%e2.15 #&#
%
\begin{eqnarray}\label{eq:2.15}
&&\cB(\a)= \bigl\{(x,y)\dvtx x\in\cR, \hg(x)+s(\cX) (x) \hz_{\a
/2}\leq y
\nonumber
\\[-8pt]
\\[-8pt]
\nonumber
&&\hspace*{75pt}\qquad
\leq\hg(x)+s(\cX) (x) \hz_{1-(\a/2)} \bigr\}.
\end{eqnarray}
The remainder of the methodology can be implemented in the following
six-step algorithm.

(1) Calculate the uncentred bootstrap residuals, $\tep_i^*=Y_i^*
-\hg
^*(X_i)$. (2) Centre them to obtain $\hep_i^*=\tep_i^*-\bep
_i^*$,
where $\bep^*=n^{-1}\sum_i \tep_i^*$. (3) Draw a double-bootstrap
resample, $\cZ^{**}=\{(X_i,Y_i^{**}),1\leq i\leq n\}$, where
$Y_i^{**}
=\hg^*(X_i)+\ep_i^{**}$ and the $\ep_i^{**}$s are sampled randomly,
with replacement, from the $\hep_i^*$s. (4) Construct the
bootstrap-world version $\cB^*(\a)$ of the band $\cB(\a)$ at
(\ref
{eq:2.15}), defined by
\[
\cB^*(\a)= \bigl\{(x,y)\dvtx x\in\cR, \hg^*(x)+s(\cX) (x) \hz _{\a/2} ^*
\leq y \leq\hg^*(x)+s(\cX) (x) \hz_{1-(\a/2)}^* \bigr\},
\]
where, reflecting (\ref{eq:2.14}), $\hz_\be^*$ is defined by
\[
P \bigl\{\hg^{**}(x)-\hg^*(x)\leq s(\cX) \hz_\be^* | \cZ^*
\bigr\}=\be,
\]
and $\cZ^*$ is defined as in step~3 of the algorithm in Section~\ref{Algorithm}. (5) For this new definition of $\cB^*(\a)$, define
$\hpi
(x,\a)$ as at (\ref{eq:2.10}). (6) Define $\hal_\xi(\a_0)$ as in
(\ref
{eq:2.12}), and take the final confidence band to be $\cB\{\hal_\xi
(\a
_0)\}$, where $\cB(\a)$ is as at (\ref{eq:2.15}).

There is also a percentile-$t$ version of this methodology, using our
quantile-based definition of $\hal_\xi(\a_0)$.

%s2.6 #&#
\subsection{\texorpdfstring{Values of $x$ that asymptotically are
covered with probability at least $1-\a_0$}{Values of x that
asymptotically are covered with probability at least
1-alphazero}}\label{Values}
Define $\|\cR\|$ to equal the Lebesgue measure of $\cR$, let $\cS$
equal the set of $x\in\cR$ such that $b(x)=0$, put $\xi_0=\|\cS\|/\|
\cR
\|$, define $\be(x,\a_0)$ to be the solution, in $\be$, of $\Phi\{
z_{1-(\be/2)}+b(x)\}-\Phi\{-z_{1-(\be/2)}+b(x)\}=1-\a_0$, and let
$\a
_\xi(\a_0)$ denote the $100 \xi\%$ quantile of values of $\be(x,\a_0)$
for $x\in\cR$. Then $\a_\xi(\a_0)$ is the solution in $\ga$ of
\[
\biggl(\int_\cR\,dx \biggr)^{ -1} \int
_\cR I \bigl\{\be(x,\a_0)\leq\ga \bigr\} \,dx=
\xi.
\]
As $\xi$ decreases, in order for the identity above to hold the value
of $\ga$ should decrease. Hence, in accordance with intuition, $\a
_\xi
(\a_0)$ decreases as $\xi$ decreases.

It can be proved that $\a_\xi(\a_0)$ is the limit in probability
of $\hal_\xi(\a_0)$. Assume that the design points $X_i$ are univariate
and that $f_X$ and $g''$ are bounded and continuous.

We showed in Section~\ref{Properties} that the naive confidence band
$\cB(\a_0)$, defined at (\ref{eq:2.2}) and having coverage $1-\a_0$,
strictly undercovers $g(x)$ when evaluated at $x$, in the asymptotic
limit, unless $b(x)=0$, and that in the latter case the coverage is
asymptotically correct, that is, equals $1-\a_0$.

Noting that $\be(x,\a_0)$ is a monotone increasing function of
$|b(x)|$, and that $b(x)=-C g''(x) f_X(x)^{1/2}$ for a positive
constant $C$, we see that if we define $\cR(\xi)=\{x\in\cR\dvtx\be
(x,\a
_0)>\a_\xi(\a_0)\}$, and $c(\xi)=\sup\{C |g''(x)| f_X(x)^{1/2}
\dvtx x\in
\cR(\xi)\}$, then the set of exceptional $x$, for which the confidence
band $\cB\{\hal_\xi(\a_0)\}$ asymptotically undercovers $(x,g(x))$, is
the set $\cS\excep$ of $x\in\cR$ such that $C |g''(x)| f_X(x)^{1/2}
>c(\xi)$. The Lebesgue measure of $\cS\excep$ equals $\max(0,\xi
-\xi
_0) \|\cR\|$. See (\ref{eq:2.2}) for a definition of $\cB(\a)$, and
step~6 of Section~\ref{Algorithm} for a definition of $\hal_\xi(\a_0)$
and a detailed account of the construction of $\cB\{\hal_\xi(\a_0)\}$.

Typically the points in $\cS\excep$ are close to peaks and troughs,
which can be identified from a graph of $\hg$. In Section~\ref{NumProps} we pay particular attention to numerical aspects of this issue.

%s2.7 #&#
\subsection{Confidence bands for probability densities}\label{ConfBands}
Analogous methods can be used effectively to construct confidence bands
for probability densities. We consider here the version of the
single-bootstrap technique introduced in Section~\ref{Algorithm}, when
it is adapted so as to construct confidence bands for densities of
$r$-variate probability distributions. Specifically, let $\cX=\{
X_1,\ldots,X_n\}$ denote a random sample drawn from a distribution with
density $f$, let $h$ be a bandwidth and $K$ a kernel, and define the
kernel estimator of $f$ by
\[
\hf(x)=\frac{1}{nh^r} \sum_{i=1}^n K
\biggl(\frac{x-X_i}{h} \biggr).
\]
This estimator is asymptotically normally distributed with
variance\break
$(nh^r)^{-1}\ka f(x)$, where $\ka=\int K^2$, and so a naive,\vadjust{\goodbreak} pointwise
confidence band for $f(x)$ is given by
%
%
%e2.16 #&#
%
\begin{eqnarray*}
&&\cB(\a)= \bigl\{(x,y)\dvtx x\in\cR, \hf(x)- \bigl[ \bigl(nh^r \bigr)
^{-1}\ka\hf(x) \bigr]^{1/2}z_{1-(\a/2)}\leq y
\\
&&\hspace*{90pt}\qquad\leq\hf(x)+ \bigl[ \bigl(nh^r \bigr)^{-1}\ka\hf(x)
\bigr]^{1/2} z_{1-(\a/2)} \bigr\};
\end{eqnarray*}
compare (\ref{eq:2.2}).

To correct $\cB(\a)$ for coverage error, draw a random sample $\cX
^*=\{
X_1^*,\ldots,\break X_n^*\}$ from the distribution with density $\hf_X$, and
define $\hf^*$ to be the corresponding kernel estimator of $\hf$,
based on $\cX$ rather than $\cX^*$,
\[
\hf^*(x)=\frac{1}{nh^r} \sum_{i=1}^n K
\biggl(\frac{x-X_i^*}{h} \biggr).
\]
Importantly, we do not generate $\cX^*$ simply by resampling
from $\cX
$. Analogously to (\ref{eq:2.9}), the bootstrap version of $\cB(\a)$
is
%
%
%e2.17 #&#
%
\begin{eqnarray*}
&&\cB^*(\a)= \bigl\{(x,y)\dvtx x\in\cR, \hf^*(x)- \bigl[ \bigl(nh^r
\bigr)^{-1}\ka\hf^*(x) \bigr]^{1/2}z_{1-(\a/2)}\leq y
\\
&&\hspace*{95pt}\qquad\leq\hf^*(x)+ \bigl[ \bigl(nh^r \bigr)^{-1}\ka\hf^*(x)
\bigr]^{1/2} z_{1-(\a
/2)} \bigr\}.
\end{eqnarray*}
For the reasons given in Remark~\ref{rem2} in Section~\ref{3Remarks} we use the
same bandwidth, $h$, for both $\cB(\a)$ and $\cB^*(\a)$.

Our bootstrap estimator $\hpi(x,\a)$ of the probability $\pi(x,\a
)=P\{
(x,f(x))\in\cB(\a)\}$ that $\cB(\a)$ covers $(x,f(x))$, is given by
$\hpi(x,\a)=P\{(x,\hg(x))\in\cB^*(\a) | \cX\}$. As in step~6
of the
algorithm in Section~\ref{Algorithm}, for a given desired coverage
level $1-\a_0$, let $\be=\hbe(x,\a_0)$ be the solution of $\hpi
(x,\be
)=1-\a_0$, and define $\hal_\xi(\a_0)$ as in (\ref{eq:2.12}). Our final
confidence band is $\cB\{\hal_\xi(\a_0)\}$. For a proportion of at
least $1-\xi$ of the values of $x\in\cR$, the limit of the probability
that this band covers $f(x)$ is not less than $1-\a_0$, and for the
remainder of values $x$ the coverage error is close to 0.

In the cases $r=1$ and 2, which are really the only cases where
confidence bands can be depicted, theoretical results analogous to
those in Section~\ref{TheoProps}, for regression, can be developed
using Hungarian approximations to empirical distribution functions.
See, for example, Theorem~3 of \citet{r34} for the case $r=1$, and
\citet{r58} and \citet{r39} for $r\geq2$. To link this argument to the
theoretical development in Appendix~B.1 in the supplementary file, we
mention that in the univariate case, the analogue of (\B.1) in that
file is
%
%
%e2.18 #&#
%
\begin{eqnarray}\label{eq:2.16}
E \bigl\{\hf^*(x) | \cZ \bigr\}-\hf(x)&=&{\tfrac{1}{2}}\ka_2
f''(x) h^2+(nh)^{-1/2}f(x)^{1/2}V(x/h)
\nonumber
\\[-8pt]
\\[-8pt]
\nonumber
&&{}+\mbox{negligible terms},
\end{eqnarray}
and (\B.3) also holds. By way of notation in (\ref{eq:2.16}) and (\B.3),
$\ka_2=\int u^2 K(u) \,du$ and, for constants $c_1$ and $c_2$, we
define $b(x)=-c_1 f''(x) f(x)^{-1/2}$ and $\De(x)=-c_2 V(x)$; and $V$
is a stationary Gaussian process with zero mean and covariance $K''*K''$.

Alternative to the definition of $\cB(\a)$ above, a confidence band
based on the square-root transform, reflecting the fact that the
asymptotic variance of $\hf$ is proportional to $f$, could be used.
Percentile and percentile-$t$ methods, using our quantile-based method
founded on $\hal_\xi(\a_0)$, can also be used.

%s3 #&#
\section{Numerical properties}\label{NumProps}

%s3.1 #&#
\subsection{Parameter settings and comparisons}\label{Parameters}

In Section~\ref{NumProps} we summarise the results of a simulation
study addressing the finite-sample performance of methodology described
in Section~\ref{Methodology}. In particular, we report empirical
coverage probabilities of nominal 95\% confidence intervals for $g(x)$,
for different~$x$, different values of $1-\xi$, different choices of
$g$, different error variances $\si^2$, and different sample sizes $n$.

For $n=100$, 200 or 400 we generated data pairs $(X_i,Y_i)$ randomly
from the model at (\ref{eq:2.1}), where the experimental errors $\ep_i$
were distributed independently as N$(0,\si^2)$ with $\si=1$, 0.5
or 0.2, and the explanatory variables $X_i$ were distributed uniformly
on $[-1,1]$. We worked with the functions $g_1$, $g_2$, and $g_3$,
defined by $g_1(x)=x+5 \phi(10 x)$, $g_2(x)=\sin(3\pi x/2) /\{1+18
x^2 (\sgn x+1)\}$ and $g_3(x)=\sin(\pi x/2) /\{1+2 x^2 (\sgn x+1)\}
$, where $\phi$ is the standard normal density and $\sgn x=1$, 0 or
$-1$ according as $x>0$, $x=0$ or $x<0$, respectively. The function
$g_1$ was used by \citet{r33}, and also by many subsequent authors;
$g_2$ is the function given by formula (7) of \citet{r2}, rescaled here
to the interval $[-1,1]$, and used extensively by \citet{r2} and in
subsequent work of other researchers; and $g_3$ is the version of $g_2$
obtained by truncating $g_2$ to the central one third of its support
interval, and rescaling so that it is supported on $[-1,1]$.

The results reported here were obtained using a standard plug-in
bandwidth, computed as suggested by \citet{r51} but employing the
variance estimator at (\ref{eq:2.8}). The cross-validation bandwidth
gives slightly better coverage results for our method, apparently
because, on average, it undersmooths a little. However, since computing
the plug-in and cross-validation bandwidths involves $O(n)$ and
$O(n^2)$ calculations, respectively, then the plug-in method is more
attractive in a numerical study that requires 1000 simulations in each
setting and sample sizes up to 400. The differences between plug-in and
cross-validation were minor in the case of competing methods since, as
discussed below, we optimised those methods over the second bandwidth.

In Section~\ref{MainResults} we report results obtained using our
method, undersmoothing without explicit bias correction, and explicit
bias correction using an oversmoothed bandwidth to estimate bias. In
the latter case we employed the regression version of a bias estimator
suggested by \citet{r54}. For each parameter setting (i.e., each sample
size $n$, each error variance $\si^2$ and each function $g_j$), when
using undersmoothing we took the bandwidth to be $\ga h$; and when
using explicit bias correction we took the bandwidth to be $h/\la$. The
values of $\ga$ and $\la$ were chosen to optimise the performance of
the two competing methods, and in particular so that those methods had
as large as possible a proportion of values $x\in\cR=[-0.9,0.9]$ that
were covered with probability at least $0.95$. To determine the best
$\ga$ and $\la$, for $n=100$ we varied $\ga$ and $\la$ in the ranges
$0.1 (0.1) 0.9$ and $0.01,0.02,0.05,0.1 (0.1) 0.9$, respectively.
For $n=200$ and 400, to reduce computation time we took the respective
ranges to be $0.2 (0.2) 1.0$ and $0.1 (0.2) 0.9$.

This approach favours the two competing methods. It is required because
there do not exist, in either case, any alternative approaches that are
even moderately widely used. Of course, this situation, which arises
because of the sheer difficulty of producing appropriate empirical
bandwidths for the competing methods, is one of the motivations for our
work. Choosing $\ga$ and $\la$ empirically, as would be necessary in
practice, would introduce significant extra variability into the
competing methodologies, and so would downgrade their performance. Even
the approach taken here, which gives competing methods every
opportunity to show their advantages, typically produces competing
techniques which perform less well than ours.

%
%f1 #&#
\begin{figure}[b]

\includegraphics{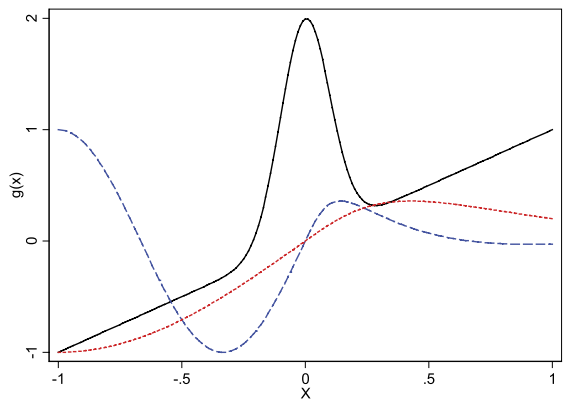}

\caption{Conditional mean functions. Solid line is $g_{1}(x)$. Long
dashes are $g_{2}(x)$. Short dashes are $g_{3}(x)$.}\label{fig:1}
\end{figure}

%s3.2 #&#
\subsection{Main results and discussion}\label{MainResults}
Graphs of $g_1$, $g_2$, and $g_3$ are shown in Figure~\ref{fig:1}. The
order $g_1,g_2,g_3$ arranges those functions in terms of decreasing
difficulty experienced by each method. In particular, $g_1$, a single
peak on a linear slope, is more challenging than $g_2$, which
represents a deep trough followed by a moderately high peak, and is
more challenging still then $g_3$, which involves a moderately steep
uphill slope followed by a gentle decrease. The extent of difficulty
can be deduced from Tables~\ref{tab:1}--\ref{tab:3}, which reveal that
the proportion of values of $x$ that are covered with probability at
least 0.95 increases, for each of the three methods, as we pass from
$g_1$ to $g_2$ and then to $g_3$.

%
%t1 #&#
\begin{table}
\caption{Simulation results for $n=100$}\label{tab:1}
\begin{tabular*}{\textwidth}{@{\extracolsep{\fill}}lcccccc@{}}
\hline
 &  &  & \multicolumn{1}{c}{$\bolds{1-\xi}$\textbf{,}} &
 \multicolumn{1}{c}{\textbf{Prop. with}} & \multicolumn{1}{c}{\textbf{Av. abs. error}} &
\multicolumn{1}{c@{}}{\textbf{Av.}} \\[-2pt]
\multicolumn{1}{@{}l}{$\bolds{\sigma}$}& \multicolumn{1}{c}{$\bolds{j}$}& \multicolumn{1}{c}{\textbf{Method}}&
\multicolumn{1}{c}{$\bolds{\gamma}$\textbf{, or} $\bolds{\lambda}$} &
\multicolumn{1}{c}{\textbf{cov. prob.} $\bolds{\geq0.95}$} & \multicolumn{1}{c}{\textbf{of cov. prob.}} & \multicolumn{1}{c@{}}{\textbf{width}} \\
\hline
1 & 1 & Ours & 0.80 & 0.685 & 0.040 & 1.172 \\
& & & 0.90 & 0.774 & 0.041 & 1.217 \\
& & & 0.95 & 0.884 & 0.042 & 1.397 \\
& 2 & & 0.80 & 0.702 & 0.025 & 0.970 \\
& & & 0.90 & 0.812 & 0.027 & 1.146 \\
& & & 0.95 & 1.000 & 0.034 & 1.322 \\
& 3 & & 0.80 & 0.945 & 0.019 & 1.009 \\
& & & 0.90 & 0.995 & 0.033 & 1.096 \\
& & & 0.95 & 1.000 & 0.042 & 1.316 \\
& 1 & Undersmooth & 0.70 & 0.801 & 0.022 & 1.105 \\
& 2 & & 0.60 & 0.840 & 0.018 & 1.076 \\
& 3 & & 0.50 & 1.000 & 0.019 & 0.989 \\
& 1 & Bias Corr. & 0.05 & 0.737 & 0.034 & 0.924 \\
& 2 & & 0.05 & 0.740 & 0.031 & 0.834 \\
& 3 & & 0.10 & 0.901 & 0.015 & 0.700 \\
0.5 & 1 & Ours & 0.80 & 0.724 & 0.038 & 0.949 \\
& & & 0.90 & 0.812 & 0.038 & 1.114 \\
& & & 0.95 & 0.895 & 0.039 & 1.197 \\
& 2 & & 0.80 & 0.823 & 0.019 & 0.822 \\
& & & 0.90 & 0.945 & 0.027 & 0.924 \\
& & & 0.95 & 0.995 & 0.034 & 0.993 \\
& 3 & & 0.80 & 0.923 & 0.018 & 0.482 \\
& & & 0.90 & 1.000 & 0.031 & 0.562 \\
& & & 0.95 & 1.000 & 0.041 & 0.642 \\
& 1 & Undersmooth & 0.80 & 0.785 & 0.024 & 0.595 \\
& 2 & & 0.70 & 0.856 & 0.018 & 0.642 \\
& 3 & & 0.70 & 1.000 & 0.019 & 0.452 \\
& 1 & Bias Corr. & 0.40 & 0.768 & 0.027 & 0.533 \\
& 2 & & 0.20 & 0.785 & 0.019 & 0.573 \\
& 3 & & 0.05 & 0.906 & 0.015 & 0.380 \\
0.2 & 1 & Ours & 0.80 & 0.409 & 0.019 & 0.421 \\
& & & 0.90 & 0.834 & 0.020 & 0.497 \\
& & & 0.95 & 0.930 & 0.027 & 0.555 \\
& 2 & & 0.80 & 0.879 & 0.020 & 0.366 \\
& & & 0.90 & 0.950 & 0.029 & 0.395 \\
& & & 0.95 & 0.961 & 0.036 & 0.424 \\
& 3 & & 0.80 & 0.945 & 0.022 & 0.231 \\
& & & 0.90 & 1.000 & 0.033 & 0.257 \\
& & & 0.95 & 1.000 & 0.041 & 0.293 \\
& 1 & Undersmooth & 0.90 & 0.801 & 0.020 & 0.399 \\
& 2 & & 0.80 & 0.818 & 0.021 & 0.282 \\
& 3 & & 0.70 & 0.978 & 0.020 & 0.217 \\
& 1 & Bias Corr. & 0.20 & 0.790 & 0.022 & 0.378 \\
& 2 & & 0.20 & 0.796 & 0.019 & 0.252 \\
& 3 & & 0.90 & 0.995 & 0.019 & 0.190 \\
\hline
\end{tabular*}
\end{table}

%
%t2 #&#
\begin{table}[t]
\caption{Simulation results for $n=200$}\label{tab:2}
\begin{tabular*}{\textwidth}{@{\extracolsep{\fill}}lcccccc@{}}
\hline
 &  &  & \multicolumn{1}{c}{$\bolds{1-\xi}$\textbf{,}} &
 \multicolumn{1}{c}{\textbf{Prop. with}} & \multicolumn{1}{c}{\textbf{Av. abs. error}} &
\multicolumn{1}{c@{}}{\textbf{Av.}} \\[-2pt]
\multicolumn{1}{@{}l}{$\bolds{\sigma}$}& \multicolumn{1}{c}{$\bolds{j}$}& \multicolumn{1}{c}{\textbf{Method}}&
\multicolumn{1}{c}{$\bolds{\gamma}$\textbf{, or} $\bolds{\lambda}$} &
\multicolumn{1}{c}{\textbf{cov. prob.} $\bolds{\geq0.95}$} & \multicolumn{1}{c}{\textbf{of cov. prob.}} & \multicolumn{1}{c@{}}{\textbf{width}} \\
\hline
1 & 1 & Ours & 0.80 & 0.745 & 0.043 & 0.967 \\
& & & 0.90 & 0.843 & 0.042 & 1.105 \\
& & & 0.95 & 0.921 & 0.043 & 1.243 \\
& 2 & & 0.80 & 0.751 & 0.023 & 0.878 \\
& & & 0.90 & 0.850 & 0.027 & 0.920 \\
& & & 0.95 & 1.000 & 0.033 & 0.962 \\
& 3 & & 0.80 & 0.900 & 0.019 & 0.734 \\
& & & 0.90 & 0.995 & 0.031 & 0.801 \\
& & & 0.95 & 1.000 & 0.041 & 0.968 \\
& 1 & Undersmooth & 0.40 & 0.989 & 0.017 & 1.266 \\
& 2 & & 0.40 & 1.000 & 0.020 & 1.228 \\
& 3 & & 0.70 & 1.000 & 0.024 & 0.545 \\
& 1 & Bias Corr. & 0.10 & 0.762 & 0.034 & 0.800 \\
& 2 & & 0.20 & 0.796 & 0.022 & 0.777 \\
& 3 & & 0.10 & 0.928 & 0.018 & 0.456 \\
\hline
\end{tabular*}
\end{table}

%
%t3 #&#
\begin{table}[b]
\caption{Simulation results for $n=400$}\label{tab:3}
\begin{tabular*}{\textwidth}{@{\extracolsep{\fill}}lcccccc@{}}
\hline
 &  &  & \multicolumn{1}{c}{$\bolds{1-\xi}$\textbf{,}} &
 \multicolumn{1}{c}{\textbf{Prop. with}} & \multicolumn{1}{c}{\textbf{Av. abs. error}} &
\multicolumn{1}{c@{}}{\textbf{Av.}} \\[-2pt]
\multicolumn{1}{@{}l}{$\bolds{\sigma}$}& \multicolumn{1}{c}{$\bolds{j}$}& \multicolumn{1}{c}{\textbf{Method}}&
\multicolumn{1}{c}{$\bolds{\gamma}$\textbf{, or} $\bolds{\lambda}$} &
\multicolumn{1}{c}{\textbf{cov. prob.} $\bolds{\geq0.95}$} & \multicolumn{1}{c}{\textbf{of cov. prob.}} & \multicolumn{1}{c@{}}{\textbf{width}} \\
\hline
1 & 1 & Ours & 0.80 & 0.746 & 0.052 & 0.963 \\
& & & 0.90 & 0.807 & 0.048 & 1.005 \\
& & & 0.95 & 0.895 & 0.046 & 1.005 \\
& 2 & & 0.80 & 0.818 & 0.022 & 0.911 \\
& & & 0.90 & 0.972 & 0.029 & 0.953 \\
& & & 0.95 & 1.000 & 0.036 & 0.953 \\
& 3 & & 0.80 & 0.840 & 0.018 & 0.907 \\
& & & 0.90 & 0.995 & 0.030 & 0.948 \\
& & & 0.95 & 1.000 & 0.041 & 0.948 \\
& 1 & Undersmooth & 0.30 & 1.000 & 0.019 & 1.208 \\
& 2 & & 0.70 & 1.000 & 0.024 & 0.637 \\
& 3 & & 0.70 & 1.000 & 0.024 & 0.429 \\
& 1 & Bias Corr. & 0.40 & 0.801 & 0.027 & 0.662 \\
& 2 & & 0.30 & 0.994 & 0.016 & 0.533 \\
& 3 & & 0.10 & 0.956 & 0.019 & 0.356 \\
\hline
\end{tabular*}
\end{table}

Table~\ref{tab:1} treats the case $n=100$, and shows, in the first
column, the values of~$\si$; in the second column, the index $j$ of the
function $g_j$; in the third column, the method; in the fourth column,
the value of $1-\xi$ (for our method), of the optimal $\ga$ (for the
undersmoothing method), and of the optimal $\la$ (for explicit bias
correction); in the fifth column, the proportion of $x\in[-0.9,0.9]$
for which the confidence band covered $g_j(x)$ with probability not
less than $0.95$ (referred to below as the ``covered proportion''); in
the sixth column, the integral average of the absolute values of
coverage errors over $x\in[-0.9,0.9]$; and in the seventh and last
column, the average widths of the confidence intervals, that is, the
average widths of the bands constructed on $\cR$. See Section~\ref{Parameters}
for definitions of $\ga$ and $\la$, and Section~\ref{Methodology} for a definition of $\xi$.

Tables~\ref{tab:2} and~\ref{tab:3} provide the same information in the
cases $n=200$ and 400, respectively, although for brevity we give
results only for $\si=1$. The numerical values in Tables~\ref{tab:1}--\ref{tab:3}
were derived by taking averages over 1000
simulations in each parameter setting. In each instance, for the sake
of brevity the tables give results only for three values of $1-\xi$,
specifically 0.8, 0.9, and 0.95. When interpreting our results, and
comparing them with those of the other methods, the reader should bear
in mind that in practice we suggest taking $1-\xi=0.9$, whereas the
competing methods have a major advantage in that we chose the tuning
parameters there to give them the largest possible value of covered proportion.

\begin{figure}
\centering
\begin{tabular}{@{}c@{}}

\includegraphics{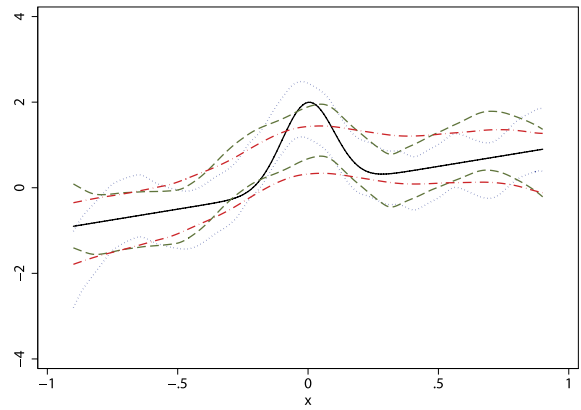}
\\
\footnotesize{(a) Proposed new method: 0.90 quantile}\\[3pt]

\includegraphics{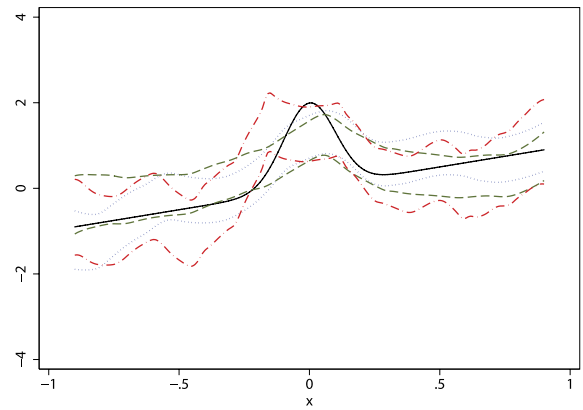}
\\
\footnotesize{(b) Conventional method with undersmoothing:
$\gamma= 0.7$}\\[3pt]

\includegraphics{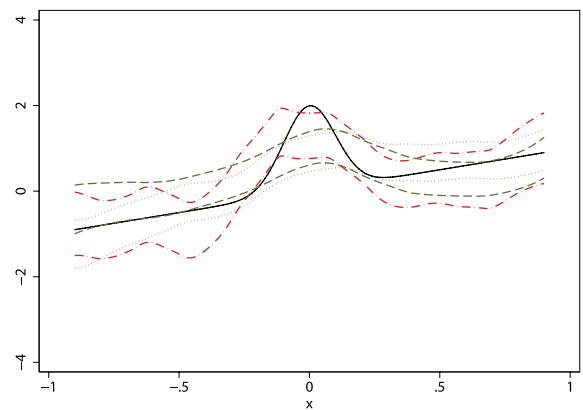}
\\
\footnotesize{(c) Conventional method with explicit bias
correction: $\lambda=0.05$}
\end{tabular}
\caption{Comparison of three methods, each panel showing three
confidence bands for interval $[-0.9,0.9]$ with $n=100$, $\sigma
^{2}=1$, and $g(x) = x + 5\phi(10x)$, $X\sim U[-1,1]$. Solid line is
$g(x)$. Lower and upper limits of the bands indicated by dashes, dots
and dash-dots.}\label{fig:2}
\end{figure}

Panels (a), (b), and (c) of Figure~\ref{fig:2} each show three typical
confidence bands in the cases of our method, of undersmoothing and of
explicit bias correction, respectively, for $g=g_1$, $n=100$ and $\si
=1$. [By ``typical'' bands we mean bands computed from the dataset for
which the integrated squared error (ISE) of the estimator took the
median value among 101 different datasets, and from the two datasets
for which ISE was closest to but not equal to the median value.] To
construct those bands in the case of our method we used $1-\xi=0.9$.
For bands in the other two cases we used the values of $\ga$ and $\la$
that maximised covered proportions in the respective parameter settings.

\begin{figure}
\centering
\begin{tabular}{@{}c@{}}

\includegraphics{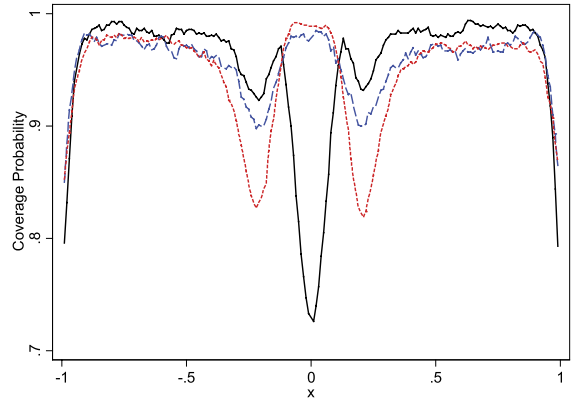}
\\
\footnotesize{(a) $g(x) = x + 5\phi(10x)$}\\[3pt]

\includegraphics{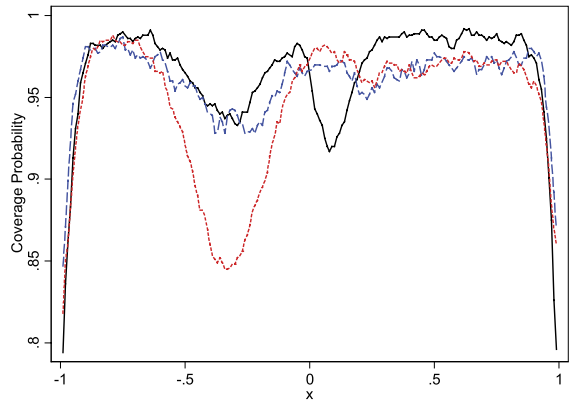}
\\
\footnotesize{(b) $g(x) = \operatorname{sin}(3\pi x/2)/\{
1+18x^{2}[\operatorname
{sgn}(x)+1]\}$}\\[3pt]

\includegraphics{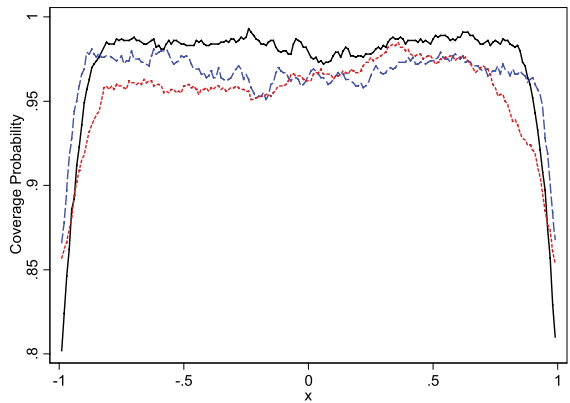}
\\
\footnotesize{(c) $g(x) = \operatorname{sin}(\pi x/2)/\{1+
2x^{2}[\operatorname
{sgn}(x)+1]\}$}
\end{tabular}
\caption{Coverage probabilities of nominal 95\% confidence band. Each
plot is for the case $n=100$, $\sigma^{2}=1$, and $X\sim U[-1,1]$, and
panels \textup{(a)}, \textup{(b)}, and \textup{(c)} are for $g = g_{1}$, $g_{2}$ and $g_{3}$,
respectively. Solid line: proposed new method. Dashes: conventional
method with undersmoothing. Dots: Conventional method with explicit
bias correction.}\label{fig:3}
\end{figure}

The three panels in Figure~\ref{fig:3} plot, as functions of $x$,
unsmoothed values of the proportions of times, out of 1000 simulations,
that the confidence band covered $(x,g(x))$. Each plot is for the case
$n=100$ and $\si=1$, and panels (a), (b), and~(c) in Figure~\ref{fig:3}
are for $g=g_1$, $g_2$ and $g_3$, respectively. The three curves in
each panel represent the method suggested in this paper, the
undersmoothing method and the explicit bias correction method,
respectively. To illustrate coverage levels at endpoints our plots
extend right across $[-1,1]$; they are not restricted to $\cR=[-0.9,0.9]$.

It can be seen from Table~\ref{tab:1} that, when $n=100$, $\si^2=1$ and
$1-\xi=0.9$, the proportion of values $x$ for which $g_j(x)$ is covered
with
%
%
%f2 #&#
%
%
%
%
%f3 #&#
probability at least 0.95, when using our method, increases from
0.77 to 0.81 and then to 0.995, for $j=1$, 2 and 3, respectively. The
corresponding values of the ``covered proportion'' are 0.80, 0.84 and 1.0
for the undersmoothing method, and 0.74, 0.74 and 0.90 in the case of
explicit bias correction. In particular, in this respect explicit bias
correction is slightly inferior to our approach, and the undersmoothing
method is slightly superior, at least in terms of the size of the
covered proportion. However, this advantage is of undersmoothing is
reversed when $\si=0.5$ or $0.2$.

In the case of undersmoothing, the value of the covered proportion can
drop sharply if there is stochastic error in choice of the bandwidth
fraction, $\ga$. Recall that in our simulation study we determine $\ga$
so that undersmoothing performs at its best, although in practice $\ga$
would be chosen implicitly using an algorithm based on estimating the
second derivative of $g_j$; this is a noisy procedure at the best of
times. To illustrate the difficulty of choosing $\ga$ in practice, we
mention that, by Table~\ref{tab:1}, when $n=100$ the optimal values of
$\ga$ are 0.7, 0.6 and 0.5 when estimating $g_1$, $g_2$ and $g_3$,
respectively, yielding covered proportions 0.801, 0.840, and 1.0,
respectively. However, if we were to mistakenly use $\ga=0.4$, 0.3 or
0.2 in these respective cases, the covered proportions would drop to
0.558, 0.354, and 0.425, respectively.

Turning to panel (b) in Figure~\ref{fig:2}, which graphs typical
confidence bands computed using the undersmoothing method, we see that
the level of undersmoothing needed to achieve a relatively high level
of covered proportion has made the band particularly wiggly, and hence
very difficult to interpret. In practice this would be quite
unsatisfactory. In comparison, the explicit bias corrected band is
about as wiggly as the band constructed using our method [compare
panels (a) and (c) in Figure~\ref{fig:2}], and both are easy to interpret.

This trend can be seen generally, for different values of $\si^2$ and
different sample sizes: The level of undersmoothing that must be used
if the undersmoothing approach is to enjoy good coverage performance,
produces bands that are distinctly unattractive because they exhibit a
high degree of spatial variability that has nothing to do with actual
features of the function $g$.

We should point out too that, in the case of undersmoothing, the
proportion of values $x\in\cR$ that are covered with probability at
least 0.95 at first increases as the bandwidth decreases, but then
starts to decrease. This is a consequence of the fact that the
confidence band quickly becomes more erratic as the bandwidth is
reduced, even more so than is shown in Figure~\ref{fig:2}. A similar
phenomenon occurs when using explicit bias correction. Here the
conservatively covered proportion of $\cR$ at first increases as we
decrease $\la$, but then it increases again. The reason is clear: If we
were to use a large bandwidth, then the bias estimator itself would be
too heavily biased, with a consequent decline in coverage performance.

The plots in Figure~\ref{fig:3} illustrate clearly the difficulty that
each approach has with the bump function $g_1$ in the interval
$(-0.3,0.3)$, where the gradient of $g_1$ changes relatively quickly.
Our approach undercovers most seriously at $x=0$, but then again, it is
honest about this; since we use $\xi=0.1$, then our approach concedes
from the outset that it can be expected to undercover approximately
10\% of points in $\cR$, and reflecting this the coverage accuracy
improves relatively quickly away from the origin. For example, it is
about 0.95 for $x=\pm0.15$, although it drops briefly down to $0.9$ in
the near vicinity of $\pm0.3$. By way of comparison, the undersmoothing
and explicit bias correction approaches perform relatively well at
$x=0$, but drop away on either side.

All three methods have less difficulty with the function $g_2$,
although it can be seen that they have more problems near the peak and
the trough than anywhere else on $\cR$. Finally, each method finds
$g_3$ relatively easy. The same trends are seen also for larger sample
sizes and smaller values of $\si$, although they are less marked in
those cases.

The average lengths of confidence bands constructed using different
methods vary in ways that are, in many instances, rather predictable.
For example, when our method produces bands with larger covered
proportion, which it does in most of the cases were considered, the
bands themselves tend to be wider, as we would expect. It is of perhaps
greater interest to focus on cases where our method has smaller covered
proportion, that is, the case $\si=1.0$ with $n=100$, 200, and 400.
When $n=100$ our bands are longer by between 7\% (in the case of $g_2$)
and 16\% (for $g_1$), despite having lower coverage. However, when
$n=200$, our bands tend to be shorter in two out of three cases (the
cases of $g_1$ and $g_2$), and when $n=400$, they are shorter in one
out of three cases (the case of $g_1$). For each method the average
lengths of bands decrease relatively slowly as sample size increases.

%s4 #&#
\section{Theoretical properties}\label{TheoProps}

%s4.1 #&#
\subsection{Theoretical background}\label{TheoBack}
In the present section we describe theoretical properties of bootstrap
methods for estimating the distribution of $\hg$. In Section~\ref{TheoPropsCE} we apply our results to underpin the arguments in
Section~\ref{Methodology} that motivated our methodology. A proof of
Theorem~\ref{thm:4.1}, below, is given in Appendix B.2 of \citet{r23}.

We take $\hg(x)$ to be a local polynomial estimator of $g(x)$, defined
by (\ref{eq:2.5}) and~(\ref{eq:2.6}). The asymptotic variance, Avar, of
the local polynomial estimator $\hg$ at $x$ is given by
%
%
%e4.1 #&#
%
\begin{equation}
\operatorname{Avar} \bigl\{\hg(x) \bigr\}=D_1 \si^2
f_X(x)^{-1} \bigl(nh_1^r
\bigr)^{-1}, \label{eq:4.1}
\end{equation}
where $D_1>0$ depends only on the kernel and $\si^2=\var(\ep)$. (If
$r=k=1$, then $D_1=\ka\equiv\int K^2$.) With this in mind we take the
estimator $s(\cX)(x)^2 \hsi^2$, introduced in Section~\ref{Properties}, of the variance of $\hg(x)$, to be $D_1 \hsi^2 \hf
_X(x)^{-1}(nh^r)^{-1}$, where $\hf_X$ is an estimator of the design
density $f_X$ and was introduced in step~1 of the algorithm in
Section~\ref{Algorithm}.

We assume that:
%
%
%e4.2 #&#
%
\begin{equation}
\label{eq:4.2} %
\begin{tabular}{p{280pt}}{(a) the data pairs
$(X_i,Y_i)$ are generated by the model at (\ref{eq:2.1}), where the design variables $X_i$ are identically
distributed, the experimental errors $\ep_i$ are identically
distributed, and the design variables and errors are totally independent; (b) $
\cR$ is a closed, nondegenerate rectangular prism in $\IR^r$; (c) the
estimator $\hg$ is constructed by fitting a local polynomial of degree $2k-1$,
where $k\geq1$; (d) $\hf_X$ is weakly and uniformly consistent, on $
\cR$, for the common density $f_X$ of the $r$-variate design
variables $X_i$; (e) $g$ has $2k$ H\" older-continuous derivatives
on an open set containing~$\cR$; (f) $f_X$ is bounded on $
\IR^r$, and H\"older continuous and bounded away from zero on an open
subset of $\IR^r$ containing $\cR$; (g) the bandwidth, $h$, used to
construct $\hg$, is a function of the data in $\cZ$ and, for constants
$C_1,C_2>0$, satisfies $P \bigl\{|h-C_1
n^{-1/(r+4k)}|>n^{-(1+C_2)/(r+4k)} \bigr\}\ra0$,\vspace*{1pt} and moreover, for constants
$0<C_3<C_4<1$, $P \bigl(n^{-C_4}\leq h\leq
n^{-C_3} \bigr)=$ $1-O \bigl(n^{-C} \bigr)$ for all $C>0$; (h)
the kernel used to construct $\hg$, at (\ref{eq:2.5}), is a spherically
symmetric, compactly supported probability density, and has $C_5$
uniformly bounded derivatives on $\IR^r$, where the positive integer
$C_5$ is sufficiently large and depends on $C_2$; and
(j) the experimental errors satisfy $E(\ep)=0$ and $E|\ep|^{C_6}<
\infty$, where $C_6>2$ is chosen sufficiently large, depending on
$C_2$.}
\end{tabular} %
\end{equation}
The model specified by (c) is standard in nonparametric regression. The
assumptions imposed in (b), on the shape of $\cR$, can be generalised
substantially and are introduced here for notational simplicity. The
restriction to polynomials of odd degree, in (c), is made so as to
eliminate the somewhat anomalous behaviour in cases where the degree is
even. See \citet{r52} for an account of this issue in multivariate
problems. Condition~(d) asks only that the design density be estimated
uniformly consistently. The assumptions imposed on $g$ and $f_X$ in (e)
and (f) are close to minimal when investigating properties of local
polynomial estimators of degree $2k-1$. Condition (g) is satisfied by
standard bandwidth choice methods, for example, those based on
cross-validation or plug-in rules. The assertion, in (g), that $h$ be
approximately equal to a constant multiple of $n^{-1/(r+2k)}$ reflects
the fact that $h$ would usually be chosen to minimise a measure of
asymptotic mean $L_p$ error, for $1\leq p<\infty$. Condition (h) can be
relaxed significantly if we have in mind a particular method for
choosing $h$. Smooth, compactly supported kernels, such as those
required by (h), are commonly used in practice. The moment condition
imposed in (j) is less restrictive than, for example, the assumption of
normality.

In addition to (\ref{eq:4.2}) we shall, on occasion, suppose that:
%
%
%e4.3 #&#
%
\begin{equation}
\label{eq:4.3} %
\begin{tabular}{p{300pt}} {the variance
estimators $ \hsi^2$ and $\hsi^*{}^2$ satisfy $P \bigl(|
\hsi- \si|>n^{-C_8} \bigr)\ra0$ and $P \bigl(|\hsi^*-\hsi|>n^{-C_8}
\bigr)\ra0$ for some $C_8>0$.}
\end{tabular} %
\end{equation}
In the case of the estimators $\hsi^2$ defined at (\ref{eq:2.7}) and
(\ref{eq:2.8}), if (\ref{eq:4.2}) holds, then so too does (\ref{eq:4.3}).

Let $h_1=C_1 n^{-1/(r+4k)}$ be the deterministic approximation to the
empirical bandwidth $h$ asserted in (\ref{eq:4.2})(g). Under (\ref
{eq:4.2}) the asymptotic bias of a local polynomial estimator $\hg$ of
$g$, evaluated at $x$, is equal to $h_1^{2k} \nabla g(x)$, where
$\nabla$ is a linear form in the differential operators $(\part/\part
x^{(1)})^{j_1}\cdots(\part/\part x^{(r)})^{j_r}$, for all choices of
$j_1,\ldots,j_r$ such that each $j_s$ is an even, positive integer,
$j_1+\cdots+j_r=2k$ [the latter being the number of derivatives assumed
of $g$ in (\ref{eq:4.2})(e)], and $x=(x^{(1)},\ldots,x^{(r)})$. For
example, if $r=k=1$, then $\nabla={\frac{1}{2}}\ka_2 (d/dx)^2$,
where $\ka
_2=\int u^2 K(u) \,du$.

Recall that $\si^2$ is the variance of the experimental error $\ep_i$.
Let $L=K*K$, denoting the convolution of $K$ with itself, and put
$M=L-K$. Let $W_1$ be a stationary Gaussian process with zero mean and
the following covariance function:
%
%
%e4.4 #&#
%
\begin{equation}
\cov \bigl\{W_1(x_1),W_1(x_2)
\bigr\}=\si^2 (M*M) (x_1-x_2).
\label{eq:4.4}
\end{equation}
Note that, since $h_1$ depends on $n$, then so too does the
distribution of $W_1$. Our first result shows that (\ref{eq:4.2}) is
sufficient for a stochastic approximation of local polynomial estimators.

%
%th1 #&#
\begin{theorem}\label{thm:4.1}
If $(\ref{eq:4.2})$ holds, then for each $n$, there exists a zero-mean
Gaussian process $W$, having the distribution of $W_1$ and defined on
the same probability space as the data $\cZ$, such that for constants
$D_2,C_7>0$,
%
%
%e4.5 #&#
%e4.6 #&#
%
\begin{eqnarray}\label{eq:4.5}
\quad&&P \Bigl[\sup_{x\in\cR} \bigl|E \bigl\{\hg^*(x) | \cZ \bigr\}-\hg(x) \nonumber
\\[-4pt]
\\[-12pt]
\nonumber
&&\hspace*{12pt}\qquad{}
-
\bigl\{h_1^{2k} \nabla g(x)
+D_2 \bigl(nh_1^r \bigr)^{-1/2}
f_X(x)^{-1/2}W(x/h_1) \bigr\}
\bigr|>h_1^{2r} n^{-C_7} \Bigr] \ra0
\end{eqnarray}
as $n\rai$. If, in addition to $(\ref{eq:4.2})$, we assume that
$(\ref
{eq:4.3})$ holds, then for some $C_7>0$,
%
%
%e4.7 #&#
%e4.8 #&#
%
\begin{eqnarray}\label{eq:4.6}
&&P \Bigl(\sup_{x\in\cR} \sup_{z\in\IR} \bigl|P \bigl[\hg
^*(x)-E \bigl\{\hg^*(x) | \cZ \bigr\}\nonumber
\\[-8pt]
\\[-8pt]
\nonumber
&&\hspace*{32pt}\qquad\leq z \bigl\{D_1 \hsi{}^2 \hf_X(x)^{-1}
\bigl(nh^r \bigr)^{-1} \bigr\}^{1/2}| \cZ \bigr] -
\Phi(z) \bigr|>n^{-C_7} \Bigr)\ra0
\end{eqnarray}
as $n\rai$.
\end{theorem}

Theorem~\ref{thm:4.1} is generically similar to other strong
approximations in the literature, although there are two differences
that are crucial to our work: the bandwidth in the theorem is a
function of the data, and has specific properties, whereas other strong
approximations in nonparametric function estimation take the bandwidth
to be deterministic; and the theorem treats data obtained using a
particular residual-based approach to resampling, and does not treat
the originally sampled data.

Result (\ref{eq:4.6}) asserts that the standard central limit theorem
for $\hg^*(x)$ applies uniformly in $x\in\cR$. In particular, the
standard deviation estimator $\{D_1 \hsi^2\*\hf_X(x)^{-1}(nh^r)^{-1}\}
^{1/2}$, used to standardise $\hg^*-E(\hg^* | \cZ)$ on the left-hand
side of~(\ref{eq:4.6}), is none other than the conventional empirical
form of the asymptotic variance of $\hg$ at (\ref{eq:4.1}), and was
used to construct the confidence bands discussed in Sections~\ref{Properties} and~\ref{Algorithm}.
The only unconventional aspect of~(\ref{eq:4.6}) is that the central limit theorem is asserted to hold
uniformly in $x\in\cR$, but this is unsurprising, given the moment
assumption in (\ref{eq:4.2})(j).

%s4.2 #&#
\subsection{Theoretical properties of coverage error}\label{TheoPropsCE}
Let $D_3=D_1^{-1/2}\si^{-1}$ and $D_4=D_2 D_3$, and define
%
%
%e4.9 #&#
%
\begin{equation}
b(x)=-D_3 f_X(x)^{1/2}\nabla g(x),\qquad
\De(x)=-D_4 W(x/h_1),\label{eq:4.7}
\end{equation}
where $W$ is as in $(\ref{eq:4.5})$. To connect these definitions to
the theoretical outline in Appendix B.1 in the supplementary file, we
note that in the present setting these are the versions of $b(x)$ and
$\De(x)$ at (\B.2) and (\B.4), respectively [$D_4 W$ in (\ref{eq:4.7})
equals $W$ in (\B.4)], and our first result in this section is a
detailed version of (\B.3):

%
%co1 #&#
\begin{corollary}\label{COR:4.1}
If $(\ref{eq:4.2})$ and $(\ref{eq:4.3})$ hold, then with $z=z_{1-(\a
/2)}$ and $b(x)$ and $\De(x)$ defined as above, we have for some $C_9>0$,
%
%
%e4.10 #&#
%
\begin{eqnarray}\label{eq:4.8}
&&P \Bigl(\sup_{x\in\cR} \bigl|\hpi(x,\a)- \bigl[\Phi \bigl\{z+b(x)+
\De(x) \bigr\} -\Phi \bigl\{-z+b(x)+\De(x) \bigr\} \bigr] \bigr|
\nonumber
\\[-4pt]
\\[-12pt]
\nonumber
&&\hspace*{245pt}\qquad>n^{-C_9}
\Bigr) \ra0
\end{eqnarray}
as $n\rai$.
\end{corollary}

Next we give notation that enables us to assert, under specific
assumptions, properties of coverage error of confidence bands. See
particularly (\ref{eq:4.13}) in Corollary~\ref{COR:4.2}, below. Results
(\ref{eq:4.11}) and (\ref{eq:4.12}) are used to derive (\ref{eq:4.13}),
and are of interest in their own right because they describe
large-sample properties of the quantities $\hbe(x,\a_0)$ and $\hal
_\xi
(\a_0)$, respectively, in terms of which our confidence bands are
defined; see Section~\ref{Algorithm}.

Given a desired coverage level $1-\a_0\in({\frac{1}{2}},1)$, define
$\hbe
(x,\a
_0)$ and $\hal_\xi(\a_0)$ as in step~6 of Section~\ref{Algorithm}, and
as at~(\ref{eq:2.12}), respectively. Let $b(x)$ and $\De(x)$
be as at~(\ref{eq:4.7}), put $d=b+\De$, and define $T=T(x,\a_0)$ to be the
solution of
\[
\Phi \bigl\{T+d(x) \bigr\}-\Phi \bigl\{-T+d(x) \bigr\}=1-\a_0.\vadjust{\goodbreak}
\]
Then $T(x,\a_0)>0$, and $A(x,\a_0)=2 [1-\Phi\{T(x,\a_0)\}]\in(0,1)$.
Define $\be=\be(x,\a_0)>0$ to be the solution of
%
%
%e4.11 #&#
%
\begin{equation}
\Phi \bigl\{z_{1-(\be/2)}+b(x) \bigr\}-\Phi \bigl\{-z_{1-(\be
/2)}+b(x)
\bigr\}=1-\a_0,\label{eq:4.9}
\end{equation}
and let $\a_\xi(\a_0)$ be the $\xi$-level quantile of the values of
$\be
(x,\a_0)$. Specifically, $\ga=\a_\xi(\a_0)$ solves the equation
%
%
%e4.12 #&#
%
\begin{equation}
\biggl(\int_\cR\,dx \biggr)^{ -1}\int
_\cR I \bigl\{\be(x,\a_0)\leq\ga \bigr\} \,dx=
\xi. \label{eq:4.10}
\end{equation}
Define $\cR_\xi(\a_0)=\{x\in\cR\dvtx I[\be(x,\a_0)>\a_\xi(\a
_0)]\}
$. Let the
confidence band $\cB(\a)$ be as at (\ref{eq:2.2}).

%
%co2 #&#
\begin{corollary}\label{COR:4.2}
If $(\ref{eq:4.2})$ and $(\ref{eq:4.3})$ hold, then for each
$C_{10},C_{11}>0$, and as $n\rai$,
%
%
%e4.13 #&#
%e4.14 #&#
%
\begin{eqnarray}
P \Bigl\{\sup_{x\in\cR\dvtx |\De(x)|\leq C_{10}} \bigl|\hbe(x,\a _0)-A(x,\a
_0)\bigr|>C_{11} \Bigr\}&\ra&0, \label{eq:4.11}
\\
P \bigl\{\hal_\xi(\a_0)\leq\a_\xi(
\a_0)+C_{11} \bigr\}&\ra&1,\label{eq:4.12}
\end{eqnarray}\vspace*{-20pt}
%
%
%e4.15 #&#
%
\begin{equation}
\begin{tabular}{p{280pt}}{for each $x\in\cR_\xi(
\a_0)$ the limit infimum of the probability $P \bigl[ \bigl(x,g(x)
\bigr)\in\cB \bigl\{\hal_\xi(\a_0) \bigr\} \bigr]$, as $n
\rai$, is not less than $1-\a_0$.} \end{tabular} %
\label{eq:4.13}
\end{equation}
\end{corollary}

Property (\ref{eq:4.12}) implies that the confidence band $\cB(\be)$,
computed using $\be=\hal_\xi(\a_0)$, is no less conservative, in an
asymptotic sense, than its counterpart when $\be=\a_\xi(\a_0)$. This
result, in company with (\ref{eq:4.13}), underpins our claims about the
conservatism of our approach. Result (\ref{eq:4.13}) asserts that the
asymptotic coverage of $(x,g(x))$ by $\cB\{\hal_\xi(\a_0)\}$ is,
for at
most a proportion $\xi$ of values of $x$, not less than $1-\a_0$.
Proofs of Corollaries~\ref{COR:4.1} and~\ref{COR:4.2} are given in
Appendix~\ref{AppA}, below.

\begin{appendix}\label{app}
%s5 #&#
\section{\texorpdfstring{Outline proofs of Corollaries \protect\ref{COR:4.1} and \protect\ref{COR:4.2}}
{Outline proofs of Corollaries 4.1 and 4.2}}\label{AppA}

%s5.1 #&#
\subsection{\texorpdfstring{Proof of Corollary \protect\ref{COR:4.1}}{Proof of Corollary 4.1}}\label{ProofCor1}

Define
\begin{eqnarray*}
\hd^*(x)&=&\frac{\hg(x)-E\{\hg^*(x) | \cZ\}}{
\{D_1 \hsi^*{}^2 \hf_X(x)^{-1}(nh^r)^{-1}\}^{1/2}},\\
 \hd(x)&=&\frac
{\hg(x)-E\{\hg^*(x) | \cZ\}}{
\{D_1 \si^2 f_X(x)^{-1}(nh_1^r)^{-1}\}^{1/2}}.
\end{eqnarray*}
Recall that, motivated by the variance formula (\ref{eq:4.1}), we take
$s(\cX)(x)^2 \hsi^2$, in the definition of the confidence band $\cB
(\a
)$ at (\ref{eq:2.2}), to be $D_1 \hsi^2 \hf_X(x)^{-1}(nh^r)^{-1}$. The
bootstrap estimator $\hpi(x,\a)$, defined at (\ref{eq:4.10}), of the
probability $\pi(x,\a)$, at (\ref{eq:2.3}), that the band $\cB(\a)$
covers the point $(x,g(x))$, is given by
%
%
%e5.1 #&#
%
\begin{eqnarray}\label{eq:A.1}
\hpi(x,\a)
&=&P \bigl\{\hg^*(x)-s(\cX) (x) \hsi^*z_{1-(\a/2)} \leq
\hg(x) \nonumber\\
&&\hspace*{22pt}{}\leq\hg^*(x)+s(\cX) (x) \hsi^*z_{1-(\a/2)} | \cZ \bigr\}
\nonumber\\
&=&P \biggl[-z_{1-(\a/2)}\leq\frac{\hg^*(x)-\hg(x)}{
\{D_1 \hsi^*{}^2 \hf_X(x)^{-1}(nh^r)^{-1}\}^{1/2}} \leq
z_{1-(\a/2)} | \cZ \biggr]
\\
&=&P \biggl[-z_{1-(\a/2)}+\hd^*(x)\leq\frac{\hg^*(x)-E\{\hg^*
(x) | \cZ\}}{
\{D_1 \hsi^*{}^2 \hf_X(x)^{-1}(nh^r)^{-1}\}^{1/2}} \nonumber\\
&& \hspace*{140pt}{}\leq z_{1-(\a/2)}+\hd^*(x) | \cZ \biggr].\nonumber
\end{eqnarray}
If both (\ref{eq:4.2}) and (\ref{eq:4.3}) hold, then by (\ref{eq:4.5}),
(\ref{eq:4.6}), (\ref{eq:A.1}), and minor additional calculations,
%
%
%e5.2 #&#
%
\begin{eqnarray}\label{eq:A.2}
&&P \Bigl(\sup_{x\in\cR} \bigl|\hpi(x,\a)- \bigl[\Phi \bigl
\{z_{1-(\a/2)}+\hd(x) \bigr\} -\Phi \bigl\{-z_{1-(\a/2)}+\hd(x) \bigr\}
\bigr]\bigr |
\nonumber
\\[-8pt]
\\[-8pt]
\nonumber
&&\hspace*{232pt}\qquad>n^{-C_9} \Bigr) \ra0.
\end{eqnarray}
Now, $-\hd(x)=D_3 f_X(x)^{1/2}\nabla g(x)+D_4 W(x/h_1)$ where
$D_3=D_1^{-1/2}\si^{-1}$ and $D_4=D_2 D_3$, and so (\ref{eq:4.8}) follows
from (\ref{eq:A.2}).

%s5.2 #&#
\subsection{\texorpdfstring{Proof of Corollary \protect\ref{COR:4.2}}{Proof of Corollary 4.2}}\label{ProofCor2}
Result (\ref{eq:4.11}) follows from (\ref{eq:4.8}). Shortly we shall
outline a proof of (\ref{eq:4.12}); at present we use (\ref{eq:4.12})
to derive (\ref{eq:4.13}). To this end, recall that $\ga=\a_\xi(\a_0)$
solves equation (\ref{eq:4.10}) when $z=z_{1-(\be/2)}$, and $\be=\be
(x,\a_0)>0$ denotes the solution of equation (\ref{eq:4.9}). If (\ref
{eq:4.12}) holds, then (\ref{eq:4.13}) will follow if we establish that
result when $\hal_\xi(\a_0)$, in the quantity $P[(x,g(x))\in\cB\{
\hal
_\xi(\a_0)\}]$ appearing in (\ref{eq:4.13}), is replaced by $\a_\xi
(\a
_0)$. Call this property (P). Now, the definition of $\a_\xi(\a_0)$,
and the following monotonicity property,
%
%
%e5.3 #&#
%
\begin{equation}
\begin{tabular}{p{280pt}} {$\Phi(z+b)-\Phi(-z+b)$ is a decreasing (resp.,
increasing) function of $b$ for $b>0$ (resp., $b<0$) and for each $z>0$,}
\end{tabular} %
\label{eq:A.3}
\end{equation}
ensure that
\[
\liminf_{n\rai} P \bigl[ \bigl(x,g(x) \bigr)\in\cB \bigl\{
\a_\xi(\a_0) \bigr\} \bigr]\geq1-\a_0
\]
whenever $\be(x,\a_0)\leq\a_\xi(\a_0)$, or equivalently, whenever
$x\in
\cR_\xi(\a_0)$. This establishes (P).

Finally we derive (\ref{eq:4.12}), for which purpose we construct a
grid of edge width~$\delta$, where $\delta$ is small [see (\ref{eq:A.4})
below], and show that if this grid is used to define $\hal_\xi(\a_0)$
[see (\ref{eq:2.12})], then (\ref{eq:4.12}) holds. Let $x_1',\ldots,x_{N_1}'$ be the centres of the cells, in a regular rectangular grid
in $\IR^r$ with edge width $\delta_1$, that are contained
within $\cR$.
(For simplicity we neglect here cells that overlap the boundaries of
$\cR$; these have negligible impact.) Within each cell that intersects
$\cR$, construct the smaller cells (referred to below as subcells) of a
subgrid with edge width $\delta=m^{-1}\delta_1$, where $m=m(\delta
_1)\geq1$ is an
integer and $m\sim\delta_1^{-c}$ for some $c>0$. Put $N=m^r N_1$; let
$x_{j\ell}$, for $j=1,\ldots,{N_1}$ and $\ell=1,\ldots,m^r$, denote the
centres of the subcells that are within the cell that has centre
$x_j'$;\vspace*{1pt} and let $x_1,\ldots,x_N$ be an enumeration of the values of
$x_{j\ell}$, with $x_{11},\ldots,x_{1m}$ listed first, followed by
$x_{21},\ldots,x_{2m}$, and so on. Recalling the definition of $\hal
_\xi
(\a_0)$ at (\ref{eq:2.12}), let $\hal_\xi(\a_0,\delta)$ denote
the $\xi
$-level quantile of the sequence $\hal(x_1,\a_0),\ldots,\hal(x_N,\a_0)$.

Let $h_1=C_1 n^{-1/(r+4k)}$ represent the asymptotic size of the
bandwidth asserted in (\ref{eq:4.2})(g), and assume that
%
%
%e5.4 #&#
%
\begin{equation}
\delta=O \bigl(n^{-B_1} \bigr),\qquad  1/(r+4k)<B_1<\infty.
\label{eq:A.4}
\end{equation}
Then
%
%
%e5.5 #&#
%
\begin{equation}
\delta=O \bigl(h_1 n^{-B_2} \bigr)\label{eq:A.5}
\end{equation}
for some $B_2>0$. In particular, $\delta$ is an order of magnitude smaller
than $h_1$.

Recall that $A(x,\a_0)=2 [1-\Phi\{Z(x,\a_0)\}]\in(0,1)$, where
$Z=Z(x,\a_0)>0$ is the solution of
\[
\Phi \bigl\{Z+b(x)+\De(x) \bigr\}-\Phi \bigl\{-Z+b(x)+\De(x) \bigr\}=1-
\a_0,
\]
and $\De(x)=-D_4 W(x/h_1)$; and that $\be=\be(x,\a_0)>0$ solves
$\Phi\{
\be+b(x)\}-\Phi\{-\be+b(x)\}=1-\a_0$. Define $e(x,\a_0)=2 [1-\Phi
\{\be
(x,\a_0)\}]$. Given a finite set $\cS$ of real numbers, let $\quant
_\xi
(\cS)$ and $\med(\cS)=\quant_{1/2}(\cS)$ denote, respectively, the
$\xi
$-level empirical quantile and the empirical median of the elements
of $\cS$. Noting~(\ref{eq:A.3}), and the fact that the stationary
process $W$ is symmetric ($W$ is a zero-mean Gaussian process the
distribution of which does not depend on $n$), it can be shown that $P\{
Z(x,\a_0)>\be(x,\a_0)\}=P\{Z(x,\a_0)\leq\be(x,\a_0)\}={\frac
{1}{2}}$. Therefore
the median value of the random variable $A(x,\a_0)$ equals $e(x,\a_0)$.
Hence, since the lattice subcell centres $x_{j1},\ldots,x_{jm^r}$ are
clustered regularly around $x_j$, it is unsurprising, and can be proved
using (\ref{eq:A.5}), that the median of $A(x_{j1},\a_0),\ldots,A(x_{jm^r},\a_0)$ is closely approximated by $e(x,\a_0)$, and in
particular that for some $B_3>0$ and all $B_4>0$,
\[
P \Bigl\{\max_{j=1,\ldots,N_1} \bigl|\med \bigl\{A(x_{j1},
\a_0),\ldots,A(x_{jm^r},\a_0) \bigr\}
-e(x_j,\a_0) \bigr| >n^{-B_3} \Bigr\}=O
\bigl(n^{-B_4} \bigr).
\]
Therefore, since the $\xi$-level quantile of the points in the set
\[
\bigcup_{j=1}^{N_1} \bigl
\{A(x_{j1},\a_0),\ldots,A(x_{jm^r},
\a_0) \bigr\}
\]
is bounded below by $\{1+o_p(1)\}$ multiplied by the $\xi$-level
quantile of the $N_1$ medians
\[
\med \bigl\{A(x_{j1},\a_0),\ldots,A(x_{jm^r},
\a_0) \bigr\},\qquad  1\leq j\leq N_1,
\]
then for all $\eta>0$,
%
%
%e5.6 #&#
%
\begin{equation}
\qquad P \bigl[\quant_{1-\xi} \bigl\{A(x,\a_0)\dvtx x\in\cR \bigr\}
\leq\quant_{1-\xi} \bigl\{e(x,\a_0)\dvtx x\in\cR \bigr\}+\eta
\bigr]\ra1.\label{eq:A.6}
\end{equation}

Since $\quant_{1-\xi} \{e(x,\a_0)\dvtx x\in\cR\}=\a_\xi(\a_0)$
then, by
(\ref{eq:A.6}),
%
%
%e5.7 #&#
%
\begin{equation}
P \bigl[\quant_{1-\xi} \bigl\{A(x,\a_0)\dvtx x\in\cR \bigr\}
\leq\a_\xi(\a_0)+\eta \bigr]\ra1.\label{eq:A.7}
\end{equation}
In view of (\ref{eq:4.11}),
%
%
%e5.8 #&#
%
\begin{equation}
\qquad P \bigl[ \bigl|\quant_{1-\xi} \bigl\{A(x,\a_0)\dvtx x\in\cR \bigr
\} -\quant_{1-\xi} \bigl\{\hbe(x,\a_0)\dvtx x\in\cR \bigr\} \bigr|>
\eta \bigr]\ra0\label{eq:A.8}
\end{equation}
for all $\eta>0$, and moreover, if $\delta$ satisfying (\ref
{eq:A.4}) is
chosen sufficiently small,
%
%
%e5.9 #&#
%
\begin{equation}
\quant_{1-\xi} \bigl\{\hbe(x,\a_0)\dvtx x\in\cR \bigr\}-
\hal_\xi(\a_0)\ra0\label{eq:A.9}
\end{equation}
in probability. [This can be deduced from the definition of $\hal_\xi
(\a
_0)$ at (\ref{eq:2.12}).] Combining (\ref{eq:A.7})--(\ref{eq:A.9}) we
deduce that $P \{\hal_\xi(\a_0)\leq\a_\xi(\a_0)+\eta\}
\ra1$ for
all $\eta>0$, which is equivalent to (\ref{eq:4.12}).
\end{appendix}

% AOS,AOAS: If there are supplements please fill:
% \sname{Supplement A}
% \stitle{Title}
% \slink[doi]{10.1214/00-AOASXXXXSUPP}
% \sdatatype{.pdf}"
% \sdescription{Some text}


\begin{supplement}[id=suppA]
\stitle{Appendix B}
\slink[doi]{10.1214/13-AOS1137SUPP} %[doi,text={...}] - jei reikia %suskaldyti doi
\sdatatype{.pdf}
\sfilename{aos1137\_supp.pdf}
\sdescription{The supplementary material in Appendix B.1 outlines
theoretical properties underpinning our methodology, while Appendix B.2
contains a proof of Theorem~\ref{thm:4.1}.}
\end{supplement}

% imsref loaded by akundreckaite, 2013-08-06 14:01:18
%

\printaddresses

\end{document}